# Analysis of Voros criterion: what derivatives involving the logarithm of the Riemann xi-function at *z=1/2* should be non-negative for the Riemann hypothesis holds true

Sergey K. Sekatskii

*Laboratoire de Physique de la Matière Vivante, IPSB, Ecole Polytechnique Fédérale de Lausanne, BSP, CH 1015 Lausanne, Switzerland*
E-mail : Serguei.Sekatski@epfl.ch

Recently, Voros has found the sums involving certain powers of $z-1/2$, which, when taken over Riemann xi-function zeroes $\rho$, must be positive for the Riemann hypothesis holds true and vice versa. Here we analyze these sums, write them as expressions involving only non-negative even powers of $\rho-1/2$, and show that the Riemann hypothesis is equivalent for the non-negativity of the derivatives $\frac{1}{(2n-1)!}\frac{d^{2n}}{dz^{2n}}(F_{2n}(z)\cdot\ln(\xi(z)))|_{z=1/2}\geq 0$ where $F_{2n}(z)=4\sum_{k=0}^{n-1}(n-k)A_{k,n}(z-1/2)^{2k}$ with the coefficients $A_{k,n}=a^{2k-2n}\sum_{l=k}^{n}C_{2n}^{2l}C_{l}^{k}$, $C_{j}^{m}$ are binomial coefficients, for any *n*=1, 2, 3… and any real *a*>1/14.

## 1. Introduction.

In 1997, Li has established the following criterion equivalent to the Riemann hypothesis concerning non-trivial zeroes of the Riemann $\zeta$-function (see e.g. [1] for standard definitions and discussion of the general properties of this function) and now bearing his name (Li's criterion) [2]:

**Li's criterion.** *Riemann hypothesis is equivalent to the non-negativity of the following numbers*

$$\lambda_n \equiv \frac{1}{(n-1)!} \frac{d^n}{dz^n}(z^{n-1} \ln(\xi(z)))|_{z=1}$$

*for any positive integer n.*

In 2013, the present author generalized Li's criterion in the following way [3]:

**Theorem 1.** *Riemann hypothesis is equivalent to the non-negativity of all derivatives $\frac{1}{(n-1)!} \frac{d^n}{dz^n}((z-a)^{n-1} \ln(\xi(z)))|_{z=1-a}$ for all positive integers n and any real a<1/2; correspondingly, it is equivalent also to the non-positivity of all derivatives $\frac{1}{(n-1)!} \frac{d^n}{dz^n}((z-a)^{n-1} \ln(\xi(z)))|_{z=1-a}$ for all positive integers n and any real a>1/2.*

Clearly, Li's criterion is a particular case *a*=0 of this Theorem.

These criteria are equivalent to the non-negativity of the sums $\Sigma_{L,a}^n = \sum_\rho (1 - \left(\frac{\rho+a}{\rho-1-a}\right)^n) \geq 0$ over all non-trivial Riemann function zeroes $\rho = \sigma_k + it_k$ taken into account their multiplicity (throughout the paper we always understand the sums over Riemann function zeroes in this sense) for all *n*=1, 2, 3… and all real *a* except *a*=1/2.

Very recently, Voros succeeded to establish inequalities involving the following sums taken over non-trivial Riemann function zeroes equivalent to the Riemann hypothesis ([4]; his original result is formulated in slightly different but equivalent form):

**Voros criterion**. *Riemann hypothesis is equivalent to the non-negativity of the following sums taken over all non-trivial Riemann function zeroes $\rho$ taking into account their multiplicity:*

$$\sum_{\rho} (1-(1+\frac{1}{2(\rho-1/2)^2}-\frac{1}{\rho-1/2}\sqrt{1+\frac{1}{4(\rho-1/2)^2}})^{-n} +$$
$$1-(1+\frac{1}{2(\rho-1/2)^2}+\frac{1}{\rho-1/2}\sqrt{1+\frac{1}{4(\rho-1/2)^2}})^{-n}) \geq 0 \ for\ n=1,\ 2,\ 3\ldots \qquad (1).$$

These sums resemble sums $\Sigma_{L,a}$ appearing above in the generalized Li's criterion, but they are taken over descending powers rather than ascending ones. Additionally, multi-valued functions (square roots) are involved here which apparently prevents a construction of derivatives corresponding to these sums in the same way as derivatives $\frac{1}{(n-1)!}\frac{d^n}{dz^n}((z-a)^{n-1}\ln(\xi(z)))|_{z=1-a}$ correspond to the sums $\Sigma_{L,a}$. (Throughout the paper we define the square root of a positive real number $a$ as a *positive* number $\sqrt{a}$. The value of the square root of an arbitrary complex number $z$ is then defined by continuation in a plane having the cut at $[0,-\infty)$).

In this Note we analyze and modify Voros' sums arriving to the following Theorem, whose proof is our main purpose:

**Theorem 2**. *Riemann hypothesis is equivalent to the non-negativity of the following sums taken over all non-trivial Riemann function zeroes $\rho$ taking into account their multiplicity:*

$$\sum_{\rho} (1-\left(\frac{1+\sqrt{1+a^2(\rho-1/2)^2}}{a(\rho-1/2)}\right)^{2n}) + (1-\left(\frac{1-\sqrt{1+a^2(\rho-1/2)^2}}{a(\rho-1/2)}\right)^{2n}) \geq 0 \qquad (2)$$

*for any real a>1/14 and any integer n=1, 2, 3… It is also equivalent to the non-negativity of the following numbers (derivatives)*

$$\frac{1}{(2n-1)!}\frac{d^{2n}}{dz^{2n}}(F_{2n}(z)\cdot \ln(\xi(z)))|_{z=1/2} \geq 0 \qquad (3)$$

*where function $F_{2n}(z)$ is a polynomial $F_{2n}(z) = 4\sum_{k=0}^{n-1}(n-k)A_{k,n}(z-1/2)^{2k}$ with the coefficients $A_{k,n} = a^{2k-2n}\sum_{l=k}^{n}C_{2n}^{2l}C_{l}^{k}$. With the same coefficients, eq. (2) can be written as*

$$\sum_{\rho}(1-(\rho-1/2)^{-2n}\sum_{k=0}^{n}A_{k,n}(\rho-1/2)^{2k}) \geq 0 \qquad (2a).$$

Here and below $C_{l}^{k} = \frac{l!}{k!(k-l)!}$ are binomial coefficients.

***Remark.*** Function $F_{2n}$ is an even function of *z-1/2*, as well as the function $\ln(\xi(z-1/2))$, which means that if we apply Leibnitz rule to the derivative $\frac{d^{2n}}{dz^{2n}}(F_{2n}(z)\cdot\ln(\xi(z)))|_{z=1/2} = \sum_{l=1}^{2n} C_{2n}^{l} \ln^{(l)}(\xi(z)) \cdot \frac{d^{2n-l}}{dz^{2n-l}}F_{2n}(z)|_{z=1/2}$, only even values of *l* matter. Indeed, any finite number of terms of the type $b_l(z-1/2)^{2l+1}$ with $b_l$=const and *l*=0, 1, 2… can be added to the function $F_{2n}$ (i.e. function $\tilde{F}_{2n}(z) = F_{2n}(z) + b_l(z-1/2)^{2l+1}$ can be used instead of $F_{2n}$ in (3)) without violating the inequality due to the evident $\frac{d^{2n}}{dz^{2n}}(b_l(z-1/2)^{2l+1}\cdot\ln(\xi(z)))|_{z=1/2} = 0$. Note also that the coefficients $A_{k,n}$ are certainly positive as well as the coefficients $4(n-k)A_{k,n}$ determining the

function $F_{2n}$. It is also clear that $\frac{d^{2k}}{dz^{2k}} F_{2n}(z)|_{z=1/2} = 4(n-k) A_{k,n}(2k)!$ for $k \le n$ and zero otherwise.

## 2. Analysis of Voros sums

André Voros noted that the mapping $f_{1a}(z) = \frac{1+\sqrt{1+a^2(z-1/2)^2}}{a(z-1/2)}$ maps the whole complex plane with a cut over a segment $\left[-it_0, it_0\right]$, $t_0$ is real $t_0 \ge 1/a$ into an exterior of the circle $|z| \le 1$ where the equality is attained on the line $\operatorname{Re} z = 1/2$, $|\operatorname{Im} z| \ge 1/a$; quite similarly, the mapping $f_{2a}(z) = \frac{1-\sqrt{1+a^2(z-1/2)^2}}{a(z-1/2)}$ maps this same plane with the cut into an interior of the circle $|z| \le 1$, and again the equality is attained on the line $\operatorname{Re} z = 1/2$, $|\operatorname{Im} z| \ge 1/a$. (Indeed, he explicitly writes expressions for a particular case $a=2$ but his paper indicates the general case quite clearly). This same observation leads Voros to establish the sums over Riemann function zeroes which need to be non-negative for the Riemann hypothesis holds true [4]. Due to the fact that there are no Riemann function zeroes with $|t| \le 14$ [1], selecting $a$ as in the conditions of Theorem 2, for all $n=1, 2, 3\ldots$, on RH all terms $1 - (\frac{1+\sqrt{1+a^2(\rho-1/2)^2}}{a(\rho-1/2)})^n$ have non-negative real part, and the same is valid for terms $1 - (\frac{1+\sqrt{1+a^2(\rho-1/2)^2}}{a(\rho-1/2)})^n$. At the same time if RH does not hold, for certain $n$ some of these terms have definitely negative (and large by module) real part.

Still we are unable to claim the non-negativity of the sums

$$\Sigma_{V,a1} = \sum_{\rho} (1-(\frac{1+\sqrt{1+a^2(\rho-1/2)^2}}{a(\rho-1/2)})^n) \text{ and}$$

$$\Sigma_{V,a}^n = \sum_{\rho} (1-(\frac{1+\sqrt{1+a^2(\rho-1/2)^2}}{a(\rho-1/2)})^n + 1 - \left(\frac{1-\sqrt{1+a^2(\rho-1/2)^2}}{a(\rho-1/2)}\right)^n) \text{ on RH, because}$$

$$\lim_{t\to-\infty} \frac{1+\sqrt{1+a^2(\rho-1/2)^2}}{a(\rho-1/2)} = -1 \quad \text{and} \quad \lim_{t\to+\infty} \frac{1-\sqrt{1+a^2(\rho-1/2)^2}}{a(\rho-1/2)} = -1 \quad \text{which}$$

means that for odd $n$ both these sums diverge. Quite the contrary, known distribution of the Riemann function zeroes $N(T) = O(T \ln T)$, where $N$ is a number of zeroes with $|\operatorname{Im}\rho| < T$ [1], guaranties the convergence of the above sums for all even $n=2, 4, 6\ldots$

This leads us to consider the sums

$$\Sigma_{V,a}^{2n} = \sum_{\rho} (1-(\frac{1+\sqrt{1+a^2(\rho-1/2)^2}}{a(\rho-1/2)})^{2n} + 1 - \left(\frac{1-\sqrt{1+a^2(\rho-1/2)^2}}{a(\rho-1/2)}\right)^{2n}) \quad (4)$$

where $n = 1, 2, 3\ldots$ We have, applying the binomial development:

$$\Sigma_{V,a}^{2n} = 2 - \frac{2}{a^{2n}(\rho-1/2)^{2n}} \sum_{l=0}^{n} C_{2n}^{2l}(1+a^2(\rho-1/2)^2)^l \text{, and now we introduce a single-}$$

valued analytical function having a pole of the $2n$-order at the point $z=1/2$:

$$f_{2n}(z) = 2 - \frac{2}{a^{2n}(z-1/2)^{2n}} \sum_{l=0}^{n} C_{2n}^{2l}(1+a^2(z-1/2)^2)^l \quad (5).$$

Using one more binomial development, we have

$$(1+a^2(z-1/2)^2)^l = \sum_{k=0}^{l} C_l^k a^{2k}(z-1/2)^{2k} \text{ and, rearranging the summation order,}$$

rewrite the function $f_{2n}(z)$ as

$$f_{2n}(z) = 2 - \frac{2}{(z-1/2)^{2n}} \sum_{k=0}^{n} A_{k,n}(z-1/2)^{2k} \quad (6)$$

where

$$A_{k,n} = a^{2k-2n} \sum_{l=k}^{n} C_{2n}^{2l} C_l^k \quad (7)$$

It might be useful to introduce $m=l-k$ and rewrite $A_{k,n} = a^{2k-2n} \sum_{m=0}^{n-k} C_{2n}^{2m+2k} C_{m+k}^k$.

Derivative of the function $f_{2n}(z)$ is equal to

$$f_{2n}'(z) = \frac{4n}{(z-1/2)^{2n+1}}\sum_{k=0}^{n}A_{k,n}(z-1/2)^{2k} - \frac{4}{(z-1/2)^{2n+1}}\sum_{k=0}^{n}A_{k,n}k(z-1/2)^{2k} =$$

$$\frac{4}{(z-1/2)^{2n+1}}\sum_{k=0}^{n-1}(n-k)A_{k,n}(z-1/2)^{2k} \qquad (8).$$

We see that $f_{2n}'(z)$ is a single-valued analytical function having an asymptotic $O(|z|^{-3})$ at large |z|. Clearly, the sums (4) are just $\sum_{\rho}f_{2n}(\rho)$ and at this stage we are able to use our standard procedure based on the generalized Littlewood theorem about contour integrals involving logarithm of an analytical function [3, 5, 6]. For completeness, we reproduce this theorem here:

**Generalized Littlewood theorem.** *Let C denotes the rectangle bounded by the lines $x = X_1$, $x = X_2$, $y = Y_1$, $y = Y_2$ where $X_1 < X_2$, $Y_1 < Y_2$ and let f(z) be analytic and non-zero on C and meromorphic inside it, let also g(z) is analytic on C and meromorphic inside it. Let F(z)=ln(f(z)), the logarithm being defined as follows: we start with a particular determination on $x = X_2$, and obtain the value at other points by continuous variation along y=const from $\ln(X_2 + iy)$. If, however, this path would cross a zero or pole of f(z), we take F(z) to be $F(z \pm i0)$ according as we approach the path from above or below. Let also the poles and zeroes of the functions f(z), g(z) do not coincide.*

*Then* $\int_C F(z)g(z)dz = 2\pi i(\sum_{\rho_g}res(g(\rho_g)\cdot F(\rho_g))) - \sum_{\rho_f^0}\int_{X_1+iY_\rho^0}^{X_\rho^0+iY_\rho^0}g(z)dz + \sum_{\rho_f^{pol}}\int_{X_1+iY_\rho^{pol}}^{X_\rho^{pol}+iY_\rho^{pol}}g(z)dz$

*where the sum is over all $\rho_g$ which are poles of the function g(z) lying inside C, all $\rho_f^0 = X_\rho^0 + iY_\rho^0$ which are zeroes of the function f(z) counted taking into account their multiplicities (that is the corresponding term is multiplied by* *m for a zero of the order m) and which lye inside C, and all $\rho_f^{pol} = X_\rho^{pol} + iY_\rho^{pol}$ which are poles of the function f(z) counted taking into account their multiplicities and which lye inside C. For this is true all relevant integrals in the right hand side of the equality should exist.*

Let us consider the same as we did before [3, 5] rectangular contour *C* with vertices at $\pm X \pm iX$ with real $X \to +\infty$, if some Riemann zero occurs on the contour just shift it a bit to avoid this. Application of the generalized Littlewood theorem to the contour integral $\int_C f_{2n}'(z)\ln(\xi(z))dz$, having in mind that $O(|z|^{-3})$ asymptotic of the function $f_{2n}'(z)$ brings its value to zero, gives $-2\pi i\sum_\rho f(\rho) + \frac{2\pi i}{(2n-1)!}\frac{d^{2n}}{dz^{2n}}(F_{2m}(z)\cdot\ln(\xi(z)))|_{z=1/2} = 0$. Here for notation purposes we have introduced a function

$$F_{2n}(z) = 4\sum_{k=0}^{n}(n-k)A_{k,n}(z-1/2)^{2k} \qquad (9).$$

Non-negativity of the sum here on RH imposes the non-negativity of the following expression: $\frac{1}{(2n-1)!}\frac{d^{2n}}{dz^{2n}}(F_{2n}(z)\cdot\ln(\xi(z)))|_{z=1/2}$ on RH, which finishes the proof of our Theorem 2.

***Remark.*** Before presenting the concluding remarks, let us demonstrate that Voros' sums (1) are exactly the sums considered here. We simply observe that

$$(\sqrt{1+\frac{1}{4(z-1/2)^2}} \pm \frac{1}{2(z-1/2)})^2 = 1 + \frac{1}{2(z-1/2)^2} \pm \frac{1}{z-1/2}\sqrt{1+\frac{1}{4(z-1/2)^2}}, \quad \text{whence}$$

$$1-(1+\frac{1}{2(z-1/2)^2} \pm \frac{1}{z-1/2}\sqrt{1+\frac{1}{4(z-1/2)^2}})^{-n} = 1-(\sqrt{1+\frac{1}{4(z-1/2)^2}} \pm \frac{1}{2(z-1/2)})^{-2n}$$

$$= 1-(\frac{\sqrt{(z-1/2)^2+1/4} \pm 1/2}{(z-1/2)})^{-2n} = 1-(\frac{1/4+(z-1/2)^2-1/4}{(z-1/2)(1/2\pm\sqrt{(z-1/2)^2+1/4})})^{-2n}$$

$$= 1-(\frac{1/2\pm\sqrt{(z-1/2)^2+1/4}}{z-1/2})^{2n}.$$

## 3. Concluding remarks

Thus we have established inequality (3) which is equivalent to the Riemann hypothesis and has a form rather similar to generalized Li's criterion (Theorem 1) but now applicable exactly at the point $z=1/2$. (Certainly, polynomial $F_{2n}$ appearing here is a bit more complicated that polynomial $(z-a)^{n-1}$ appearing in the conditions of Theorem 1, but this does not seem too principal).

Note, that in 1999 Pustyl'nikov argued that all even derivatives of the Riemann $\xi$-function (not of the $\ln\xi$-function) taken at $z=1/2$ should be positive for the Riemann hypothesis holds true and showed that they are indeed positive (unfortunately, this is not enough to prove the Riemann hypothesis) [7]; similar theorem has been proven by Coffey in 2004 [8]. Of course, connection between the derivatives of $\xi$- and $\ln\xi$-functions is complicated, but the present author has a strong feeling that this line of researches might be truly perspective one.

Another question which appears here is that of an *arithmetic interpretation* of inequality (3), similarly as this has been done for Li's [9] and generalized Li's criteria [3]. Certainly, well known property [1]

$$\ln\varsigma(z) = -\sum_p \ln\left(1 - \frac{1}{p^z}\right) = \sum_p \sum_{n=1}^{\infty} \frac{1}{np^{nz}} = \sum_{n=2}^{\infty} \frac{\Lambda(n)}{\ln n \cdot n^z} \qquad (10)$$

where Re$z>1$ and we have a sum over primes or use the van Mandgoldt function, lies in the very core here, and we already argued [10] that such an interpretation can be established not only following Weil's analysis of the explicit formulae in the number theory [11] (see also [12]), as this was done in [3, 9], but also by substituting expansion (10) and formula $\xi(z) = \frac{1}{2}z(z-1)\pi^{-z/2}\Gamma(z/2)\varsigma(z)$ into the corresponding relations between the sums and derivatives involved. In particular, in such a way we have demonstrated that if we assume that there exists a real number $1/2 < \sigma_0 < 1$, such that for all Riemann function zeroes one has $\text{Re}\,\rho \le \sigma_0 + \varepsilon$ (where $\varepsilon$ is an arbitrary small fixed positive number), the following limit holds: $\frac{\varsigma'(a)}{\varsigma(a)} = -\lim_{N \to \infty} (\sum_{m \le N} \frac{\Lambda(m)}{m^a} - \frac{N^{1-a}}{1-a})$ [13]. Differentiation of this equality with respect to $a$ readily gives a number of equalities involving higher order derivatives $\frac{d^j}{ds^j} \frac{\varsigma'(s)}{\varsigma(s)}|_{s=a}$ and sums $\sum_{m \le N} \frac{\Lambda(m) \ln^j m}{m^a}$, which can be substituted into above relations which express the sums over Riemann function zeroes via certain derivatives involving the logarithm of the Riemann zeta-function. The question either analogous relation, viz. $\frac{\varsigma'(1/2)}{\varsigma(1/2)} = -\lim_{N \to \infty} (\sum_{m \le N} \frac{\Lambda(m)}{\sqrt{m}} - 2\sqrt{N})$, on RH holds for $a=1/2$, remains unclear for the present author (probably not), so we will not sought for the arithmetic interpretation in this Note.

Finally, we could note that, similarly to Bombieri-Lagarias generalization of Li's criterion for the case of arbitrary complex number multisets [9], the above criterion also can be generalized along the same lines, sf. [3]. For this we need to consider slightly more general sums

$$\Sigma = \sum_{\rho} (1 - (\frac{1 + \sqrt{1 + a^2(\rho - \sigma)^2}}{a(\rho - \sigma)})^{2n})$$ or

$$\Sigma = \sum_{\rho} (1 - (\frac{1 + \sqrt{1 + a^2(\rho - \sigma)^2}}{a(\rho - \sigma)})^{2n} + 1 - (\frac{1 - \sqrt{1 + a^2(\rho - \sigma)^2}}{a(\rho - \sigma)})^{2n})$$, where $\sigma$ is real

and we search to establish that for all members of the multiset, composed by numbers $\rho$, $\text{Re}\,\rho = \sigma$. We also need to suppose that for the multiset at question there exists the smallest value of $|t_k| = \varepsilon > 0$ (so that we can find a suitable $a \geq 1/\varepsilon$) as well as the convergence of appropriate sums. We will not pursue this line of researches here as well as different possible generalizations of the present approach for other zeta-functions.

ACKNOWLEDGEMENT

The author thanks André Voros for information about paper [4] and useful discussions.